\documentclass{article}
\usepackage{amsmath,amsfonts,amsthm,epsfig,amssymb,color}
\newtheorem{theorem}{Theorem}
\newtheorem{lemma}[theorem]{Lemma}

\newtheorem*{claim*}{Claim}
\newtheorem{corollary}[theorem]{Corollary}

\renewcommand\Pr{{\mathop{\mathbb P{}}\nolimits}}
\renewcommand\phi{\varphi}
\newcommand\eps{\varepsilon}

\newcommand\tg{{\tilde g}}
\newcommand\dd{{\mathrm d}}
\newcommand\pmax{p_{\max}}
\newcommand\hmax{h_{\max}}

\newcommand\E{\operatorname{\mathbb E{}}}

\newcommand\norm[1]{||#1||}

\newcommand\Ec{E_3^\mathrm{crude}}

\begin{document}
\title{Erratum: Percolation on random Johnson--Mehl tessellations and related models}

\date{May 8, 2009}

\author{B\'ela Bollob\'as%
\thanks{Department of Pure Mathematics and Mathematical Statistics,
Wilberforce Road, Cambridge CB3 0WB, UK}
\thanks{Department of Mathematical Sciences,
University of Memphis, Memphis TN 38152, USA}
\thanks{Research supported in part by NSF grants DMS-0505550,
 CNS-0721983 and CCF-0728928, and ARO grant W911NF-06-1-0076}
\and Oliver Riordan%
\thanks{Mathematical Institute, University of Oxford, 24--29 St Giles', Oxford OX1 3LB, UK}
}

\maketitle

The proof presented in~\cite{JM} of the result that the critical
probability for percolation on a random Johnson--Mehl tessellation is
$1/2$ contains a (glaring!) error; we are very grateful to Rob van den Berg
for bringing this to our attention. Fortunately, the error is easy to correct;
as is often the case when one applies sharp-threshold results
such as those of Talagrand~\cite{Tal94} or Friedgut and Kalai~\cite{FK}
(in both cases extending results of Kahn, Kalai and Katznelson~\cite{KKL}
and Bourgain, Kahn, Kalai, Katznelson and Linial~\cite{BKKKL})
to `symmetric' events,
to obtain a sharp threshold one needs only enough symmetry to ensure
that many variables are equivalent, rather than total symmetry.

Let $\Pr_{p_-,p_+}^n$ denote the probability measure on $\{-1,0,1\}^n$ in
which each coordinate is independent, and is equal to $+1$ with probability $p_+$, 
and to $-1$ with probability $p_-$. An event $E\subset \{-1,0,1\}^n$ is {\em increasing}
if whenever $x\in E$ and $x\le x'$ holds coordinatewise, then $x'\in E$.
We say that $E$ has {\em symmetry of order $m$} if there is a group
action on $[n]=\{1,2,\ldots,n\}$ in which each orbit has size at
least $m$, such that the induced action on $\{-1,0,1\}^n$ preserves $E$. To
correct the proof in~\cite{JM} we need the following lemma.

\begin{lemma}\label{l1}
There is an absolute constant $c_3$ such that if $0<q_-<p_-<1/e$, $0<p_+<q_+<1/e$,
$E\subset \{-1,0,1\}^n$ is increasing and has symmetry of order $m$, and
$\Pr_{p_-,p_+}^n(E)>\eta$, then $\Pr_{q_-,q_+}^n(E)>1-\eta$ whenever
\begin{equation}\label{lb}
 \min\{q_+-p_+,p_--q_-\} \ge c_3\log(1/\eta) \pmax\log(1/\pmax) /\log m,
\end{equation}
where $\pmax=\max\{q_+,p_-\}$.
\end{lemma}
Using this result in place of Theorem 2.2 of~\cite{Voronoi} (which is simply
the special case when $m=n$, i.e., when $E$ is symmetric), the proof in~\cite{JM}
may be corrected with essentially no changes.  Indeed, the
event $E_3$ considered at the bottom of page 329, that
some $3s/4$ by $s/12$ rectangle in ${\mathbb T}(s)$ has a robustly
black horizontal crossing, is symmetric under translations of 
the space ${\mathbb T}(s)\times [0,s]$ in which the Poisson point processes
live through
vectors of the form $(x,y,0)$. Hence the corresponding discrete event
$\Ec$ considered on the next page has symmetry of order
$m=(s/\delta)^2$. To deduce Theorem 6 of~\cite{JM} from Lemma~\ref{l1}
one needs the inequality in the middle of page 330 of~\cite{JM}, but with $\log N$
replaced by $\log m$. Since $N=(s/\delta)^3$, this corresponds
simply to a change in the constant, and all remaining calculations
are unaffected.

To prove Lemma~\ref{l1} one needs a suitable influence result. Such a
result was proved for a product of $2$-element spaces by
Talagrand~\cite{Tal94}; later, Friedgut and Kalai~\cite{FK}
used a different method to obtain slightly weaker results.
One can adapt Talagrand's proof to the $3$-element
setting, obtaining a slightly weaker form of Lemma~\ref{l1} (see
the remark at the end of this note), but it seems easier to
follow the method of~\cite{FK}. Unfortunately, 
even in the two element case, Friedgut and Kalai
did not prove quite the result we need,
although their method gives it.

Given a function $f$ on a product probability space $\Omega^n$,
let $I_f(k)$ denote the {\em influence} of the $k$th coordinate
with respect to $f$, i.e., the probability of the set
of configurations $\omega$ with the property that there is some
$\omega'$ differing from $\omega$
only in the $k$th coordinate for which
$f(\omega')\ne f(\omega)$. For $A\subset \Omega^n$,
let $I_A(k)=I_f(k)$ where $f$ is the characteristic function of $A$.

Following the notation of Friedgut and Kalai~\cite{FK}, let
$V_n(p)$ denote the {\em weighted cube}, that is the $n$th power
of the $2$-element probability space in which $\Pr(0)=1-p$ and $\Pr(1)=p$.
Bourgain, Kahn, Kalai, Katznelson and Linial~\cite{BKKKL} showed that
if $f$ is any $0/1$-valued function on the $n$th power of a probability space,
then some influence $I_f(k)$ is at least a constant times $t\log n/n$,
where $t=\min\{\Pr(f^{-1}(0)),\Pr(f^{-1}(1))\}$.
Friedgut and Kalai~\cite{FK} adapted their proof to prove two extensions
(Theorems 3.1 and Theorem 3.4 in~\cite{FK}). The following result
combines these extensions.  It is also implied by Corollary 1.2
of Talagrand~\cite{Tal94}; see the remark below.

\begin{lemma}\label{l2el}
Let $0<p\le 1/2$ and let $f:V_n(p)\to \{0,1\}$ with $\Pr(f^{-1}(1))=t$.
If $I_f(k)\le \delta$ for every $k$ then $I_f=\sum_{k=1}^n I_f(k)$
satisfies the inequality
\[
 I_f \ge c\frac{t(1-t)}{p\log(1/p)} \log\left(\frac{ct(1-t)}{\delta^{1/2}I_f}\right),
\]
where $c>0$ is an absolute constant. In particular, if for some $a\le 1/2$
we have $I_f(k)\le a p^2\log(1/p)^2$ for every $k$, then
\[
 I_f \ge c'\frac{t(1-t)}{p\log(1/p)} \log(1/a),
\]
where $c'>0$ is an absolute constant.
\end{lemma}
In the related results in~\cite{FK}, it is assumed that $t\le 1/2$, in which
case the factor $(1-t)$ can be dropped. Of course this makes no
difference; indeed, $t(1-t)$ can be replaced by $\min\{t,(1-t)\}$
above, changing the constants appropriately.
The condition $p\le 1/2$
can be replaced by $p\le 1-\eps$ for any constant
$\eps>0$, but in any case the main interest is when $p$ is small.
(Such a condition is assumed implicitly in~\cite{FK}.)

Although Lemma~\ref{l2el} is not given in~\cite{FK}, it might as well have been:
to prove it one simply combines the two modifications to the
Bourgain, Kahn, Kalai, Katznelson and Linial~\cite{BKKKL} proof
that Friedgut and Kalai~\cite{FK} gave; these modifications can be applied
simultaneously without any problems.

\begin{proof}
As in~\cite{FK}, we phrase the proof in terms of modifications to that
in~\cite{BKKKL}; what follows is not intended to be read on its own.

As in~\cite{BKKKL}, the first step is to replace each factor $V_1(p)$
in the product space $V_n(p)$
by the probability space $Y=\{0,1\}^m$ with uniform measure; one can
assume that $p$ is a dyadic rational, choose $m$ so that $2^mp$ is an integer,
and take the first $(1-p)2^m$ points of $Y$ (in the binary order) to correspond
to $0\in V_1(p)$ and the last $p2^m$ to $1\in V_1(p)$.
Then, as noted in~\cite{FK}, for any function $f:V_1(p)\to \{0,1\}$,
the sum $w(f)$ of the influences of the corresponding function on $Y$
satisfies
\begin{equation}\label{wf}
 w(f)\le c_1 p\log(1/p)
\end{equation}
for some absolute constant $c_1$.
Using this in place of the bound $w(f)\le 2$ one can replace
relation (14) of~\cite{BKKKL} by $||W_k||_2^2\le c_1 p\log(1/p)I_f(k)$.

Writing $\delta_k$ for $I_f(k)$ and
using the first part of (15) in~\cite{BKKKL},
it then follows that more than half the weight of the sum
\[
 t(1-t) = ||f-\E f||_2^2= \sum_{S_1\subset [m],\dots,S_n\subset [m]: |S_1|+\cdots+|S_n|>0} {\hat f}^2(S_1,\ldots,S_n)
\]
is concentrated on terms ${\hat f}^2(S_1,\ldots,S_n)$
with
\begin{equation}\label{e1}
 0< \sum |S_i| \le 3c_1(t(1-t))^{-1} p\log(1/p)\sum_k \delta_k.
\end{equation}
[We have slightly modified the argument in~\cite{FK} to exclude the term
with all $S_i$ empty. This seems to be needed later to correct
a trivial error in~\cite{FK}.]

On the other hand, as noted by Friedgut and Kalai~\cite{FK}, with $\eps=1/\sqrt{3}$
relations (18) and (19) in~\cite{BKKKL} give
\[ 
 \sum_{k=1}^n  ||T_\eps R_k||_2^2 \le \sum_{k=1}^n ||R_k||_{4/3}^2 \le
 \sum_{k=1}^n (3\delta_k)^{3/2}, 
\]
and it follows from \cite[(20)]{BKKKL} that more than half the weight
of $||f-\E f||_2^2$ is on terms with
\begin{equation}\label{e2}
 (|S_1|+\cdots+|S_n|) \eps^{2|S_1|+\cdots+2|S_n|} \le c_2^{-1} (t(1-t))^{-1} \sum_{k=1}^n \delta_k^{3/2},
\end{equation}
for some absolute constant $c_2>0$.
Hence, some weight of $||f-\E f||_2^2$ sits where both \eqref{e1} and \eqref{e2} hold,
so these
inequalities can hold simultaneously.

Let $s$ denote a value of $|S_1|+\cdots+|S_n|$ for which
both \eqref{e1} and \eqref{e2} hold.
Since $s\ge 1$ we have $s\eps^{2s}=s3^{-s}\ge 3^{-s}$, so from \eqref{e2}
\[
 3^{-s} \le c_2^{-1} (t(1-t))^{-1} \sum_{k=1}^n \delta_k^{3/2},
\]
i.e., 
\[
 s\ge \log\left(c_2t(1-t)/\sum \delta_k^{3/2}\right)/\log 3.
\]
Combined with \eqref{e1}, this gives
\[
  3c_1(t(1-t))^{-1} p\log(1/p)\sum_k \delta_k \ge \log\left(c_2t(1-t)/\sum \delta_k^{3/2}\right)/\log 3,
\]
i.e.,
\begin{equation}\label{final}
 I_f \ge c_3 \frac{t(1-t)}{p\log(1/p)} \log\left(c_2t(1-t)/\sum \delta_k^{3/2}\right)
\end{equation}
for some absolute $c_3>0$, where $I_f=\sum \delta_k$ is the sum of the influences.

Note that \eqref{final} is valid for any $0/1$-valued function on $V_n(p)$.
Assuming now that $\delta_k\le \delta$ for all $k$,
we have
$\sum \delta_k^{3/2}\le \delta^{1/2}\sum\delta_k=\delta^{1/2}I_f$, so
\[
 I_f \ge c_3 \frac{t(1-t)}{p\log(1/p)} \log\left(\frac{c_2t(1-t)}{\delta^{1/2}I_f}\right),
\]
proving the first part of the result.

Define $a$ by $\delta= a (p\log(1/p))^2$, and set $x=t(1-t)/(p\log(1/p))$.
Suppose that $I_f\le b x$.
Then we have
\[
 I_f \ge c_3 x \log\left(\frac{c_2t(1-t)}{a^{1/2}p\log(1/p) b t(1-t)(p\log(1/p))^{-1}}\right)
 = c_3 x \log\left(\frac{c_2}{a^{1/2}b}\right).
\]
Since $I_f=bx$ it follows that $b \ge c_3\log(c_2/(a^{1/2}b))$.
Assuming that $a\le 1/2$, say,
it follows that 
\[
 I_f\ge 2c_4 \log(a^{-1/2}) x =  c_4 \frac{t(1-t)}{p\log(1/p)} \log(1/a)
\]
for some absolute constant $c_4>0$, as claimed.
\end{proof}

As noted in~\cite{Voronoi}, because of the form of the proof, the extension
to the probability space $W_{p_-,p_+}^n$, i.e., $\{-1,0,1\}^n$ with
the product measure $\Pr_{p_-,p_+}^n$, is immediate.
We state only the second part.

\begin{corollary}\label{c}
For every $0<p_-,p_+\le 1/e$ and
every function $f:W_{p_-,p_+}^n\to \{0,1\}$ with $\Pr(f^{-1}(1))=t$,
if $a\le 1/2$ and $I_f(k)\le a \pmax^2\log(1/\pmax)^2$ for every $k$, then
\[
 I_f \ge c\frac{t(1-t)}{\pmax\log(1/\pmax)} \log(1/a),
\]
where $\pmax=\max\{p_-,p_+\}$ and $c>0$ is an absolute constant.
\end{corollary}
\begin{proof}
The proof is almost identical to that of Lemma~\ref{l2el};
the first step is to replace each factor $W_{p_-,p_+}$ by
$Y=\{0,1\}^m$, noting that this time one has
\[
 w(f)\le cp_+\log(1/p_+) +cp_-\log(1/p_-) \le 2c\pmax\log(1/\pmax)
\]
in place of \eqref{wf}. From this point on the original probability space
is irrelevant.
\end{proof}

Using standard methods, it is easy to deduce Lemma~\ref{l1}
from Corollary~\ref{c}.
\begin{proof}[Proof of Lemma~\ref{l1}]
Since the left-hand side of \eqref{lb} is at most $\pmax$,
taking $c_3$ large we may assume that $\log m\ge 100\log(1/\pmax)$, say,
i.e., that $m\ge \pmax^{-100}$.

For $0\le h\le \hmax=\min\{q_+-p_+,p_--q_-\}$
let $r_+=p_++h$ and $r_-=p_--h$, let $g(h)=\Pr_{r_-,r_+}^n(E)$,
and let $\tg(h)=\log(g(h)/(1-g(h)))$.
Note that $g(0)=\Pr_{p_-,p_+}^n(E)\ge \eta$, so $\tg(0)\ge -\log(1/\eta)$,
while $\Pr_{q_-,q_+}^n(E)\ge g(\hmax)$.
We claim that
\begin{equation}\label{claim}
 \frac{\dd\tg}{\dd h} \ge \frac{2\log m}{c_3\pmax\log(1/\pmax)}
\end{equation} for any $0\le h\le \hmax$.
Assuming this, using the lower bound \eqref{lb} on $\hmax$,
we then have $\tg(\hmax)\ge \tg(0)+2\log(1/\eta)\ge \log(1/\eta)$,
giving $g(\hmax)>1-\eta$, and hence $\Pr_{q_-,q_+}^n(E)>1-\eta$, as
required.

To prove \eqref{claim}, note that $\dd\tg/\dd h = (g(1-g))^{-1}\dd g/\dd h$,
and that, by a form of the Margulis--Russo formula,
the derivative of $g(h)$ is at least $I_f=\sum_k I_f(k)$, where
$f$ is the characteristic function of $E$ and we evaluate
the influences in the product space $\Pr_{r_-,r_+}^n$.
Hence \eqref{claim} follows if we can show that
\begin{equation}\label{need}
 I_f \ge 2c_3^{-1} \frac{t(1-t)}{\pmax\log(1/\pmax)}\log m,
\end{equation}
where $t=g(h)=\Pr_{r_-,r_+}^n(E)$.

Suppose first that some influence $I_f(k)$ is at least $m^{-1/2}$, say.
Then, from the symmetry assumption, at least $m$ influences
are at least this large, and $I_f\ge m^{1/2}\ge m^{1/3}\pmax^{-2}$.
Taking $c_3$ large enough, this is much larger than the bound
in \eqref{need}. (The factor $t(1-t)$ works in our favour.)
On the other hand, if $I_f(k)\le m^{-1/2}$ for all $k$,
then $a=\max I_f(k)\pmax^{-2}\log(1/\pmax)^{-2}\le m^{-1/3}$, say,
and Corollary~\ref{c} gives \eqref{need}.
\end{proof}

Let us remark briefly on the relationship of the results
above to those of Talagrand~\cite{Tal94}. Note first
that given an increasing subset $A$ of
the weighted cube $V_n(p)$, the quantity
$\mu_p(A_i)$ in~\cite{Tal94}
is exactly $pI_A(i)$, where $I_A(i)=I_f(i)$ with
$f$ the characteristic function of $A$.
Theorem 1.1 of~\cite{Tal94} thus states in the notation above that for
some universal constant $c>0$ and any $A\subset V_n(p)$
with $\Pr(A)=t$,
\[
 \sum_{k=1}^n \frac{p(1-p)I_A(k)}{\log[1/(p(1-p)I_A(k))]} \ge c \frac{t(1-t)}{\log[2/(p(1-p))]},
\]
and Corollary 1.2 in~\cite{Tal94} gives
\begin{equation}\label{T2}
 I_A \ge \frac{ct(1-t)}{p(1-p)\log[2/(p(1-p))]} \log(1/\eps)
\end{equation}
whenever $p(1-p)I_A(i)\le \eps$ for all $i$. (In fact, we have reinserted
an irrelevant factor $(1-p)$ omitted in~\cite{Tal94}. Also, Talagrand
states his results for monotone subsets, but this condition is not used.)
These results immediately imply Lemma~\ref{l2el}; indeed, Theorem 1.1 of~\cite{Tal94}
is stronger.
The corollary \eqref{T2} is superficially stronger than Lemma~\ref{l2el}, but in practice
most likely exactly equivalent, as in the applications one always assumes that $\eps$
is smaller than some large power of $p$, and then the apparent differences
are irrelevant up to changing the constants.

Why then did we start from (a form of) the Friedgut--Kalai result instead of Talagrand's?
The answer is that the extension to a power of a 3-element space is clearer, at least to us.
One approach is as follows. First, consider any function
$f$ defined on the {\em weighted cube} $V_n(p_1,\ldots,p_n)$,
i.e., the product of the probability spaces $V_1(p_1),\ldots,V_1(p_n)$, with state space $\{0,1\}^n$.
The proof of Lemma~\ref{l2el} gives the following result; note that
assuming $p_i\le 1/2$ loses no generality, as one can replace $p_i$ by $1-p_i$.
\begin{lemma}\label{l1dif}
Let $0<p_1,\ldots,p_n\le 1/2$ and let
$f:V_n(p_1,\ldots,p_n)\to \{0,1\}$ with $\Pr(f^{-1}(1))=t$.
If $I_f(k)\le \delta$ for every $k$ then $I_f=\sum_{k=1}^n I_f(k)$
satisfies the inequality
\[
 I_f \ge c\frac{t(1-t)}{\max_i p_i\log(1/p_i)} \log\left(\frac{ct(1-t)}{\delta^{1/2}I_f}\right),
\]
where $c>0$ is an absolute constant.
\end{lemma}

Turning to Talagrand's version,
for $x\in \{0,1\}^n$, write $U_i(x)$ for the point obtained by changing the $i$th coordinate
of $x$, and, adapting Talagrand's definition in the obvious way,
set $\Delta_i f(x)=(1-p_i)(f(x)-f(U_i(x)))$ if $x_i=1$ and
$\Delta_i f(x)=p_i(f(x)-f(U_i(x)))$ if $x_i=0$.
The proof of Theorem 1.5 in~\cite{Tal94} goes through {\em mutatis mutandis}
to give the following result.
\begin{theorem}
There is an absolute constant $K$ such that, for any function $f$ on any
space $V_n(p_1,\ldots,p_n)$ with $\E f=0$, we have
\[
 \norm{f}_2^2 \le K \log\left(\frac{2}{\min_i p_i(1-p_i)}\right)
  \sum_i \frac{ \norm{\Delta_i f}_2^2} { \log(e\norm{\Delta_i f}_2/\norm{\Delta_i f}_1) },
\]
where the expectation $\E$ and norms $\norm{\cdot}_q$ are calculated with respect
to the probability measure on $V_n(p_1,\ldots,p_n)$.
\end{theorem}
Since the modifications to Talagrand's proof are essentially trivial, we omit the
details. The key point is that, much of the time, one coordinate is considered
at a time, and one should simply replace $p$ by the corresponding $p_i$
wherever it appears (for example, in the definition of $r_i(x)$).
The minimum over $i$ comes in when the upper bound $|r_i(x)-r_i(y)|\le \theta$
on page 1580 of~\cite{Tal94} is used;
this now requires $\theta=\max_i 1/\sqrt{p_i(1-p_i)}$.

Translating back to influences in the case where $f$ is the characteristic
function of $A$ with the constant $\Pr(A)$ subtracted, one has
\[
\norm{\Delta_i f}_q^q = I_A(i)(p_i(1-p_i)^q +(1-p_i)p_i^q),
\]
giving
$\norm{\Delta_i f}_2 = \sqrt{I_A(i)p_i(1-p_i)}$ and
$\norm{\Delta_i f}_1 = 2I_A(i)p_i(1-p_i)$.
With $\eps=\max_i p_i(1-p_i)I_A(i)\le \max I_A(i)$, the equivalent of Corollary 1.2
in~\cite{Tal94} one obtains is thus that
\[
 I_A \ge \frac{ct(1-t)} { (\max_i p_i(1-p_i)) \log[ 2/(\min_i p_i(1-p_i))] } \log(1/\eps).
\]
If the $p_i$ vary wildly, this inequality
seems to be weaker then Lemma~\ref{l1dif}, although the
difference may not matter.  If the maximum and minimum are close, it is
stronger than Lemma~\ref{l1dif}, although most likely equivalent
for almost all applications.

To apply this result to~\cite{JM},
roughly speaking one can replace each copy of $\{-1,0,1\}$
by $\{0,1\}^2$ with an appropriate measure, with $p_1$ close to $1-p_-$ and $p_2$ close
to $p_+$. Since in this case $p_-$ and $p_+$ are close, while we do not quite
obtain Lemma~\ref{l1}, we obtain a result that is good enough for our application.
Since we have given a different proof above, we omit the details.

\end{document}